\documentclass[12pt]{article}
\usepackage{amsfonts}

\usepackage{latexsym}
\usepackage{amssymb}
\usepackage{epsfig}
\usepackage{amsmath}
\usepackage{amsthm}
\usepackage[mathscr]{eucal}
\usepackage{graphicx}

\usepackage{caption}
\usepackage{subcaption}

\newtheorem{Lemma}{Lemma}

\newtheorem{Theorem}[Lemma]{Theorem}

\newtheorem{Definition}{Definition}

\renewcommand{\qed}{\hfill{\ \ \rule{2mm}{2mm}} \vspace{0.2in}}

\begin{document}

\title{Duality in percolation via outermost boundaries II: Star connected components and left right crossings}
\author{ \textbf{Ghurumuruhan Ganesan}
\thanks{E-Mail: \texttt{gganesan82@gmail.com} } \\
\ \\
New York University, Abu Dhabi}
\date{}
\maketitle

\begin{abstract}



Tile \(\mathbb{R}^2\) into disjoint unit squares \(\{S_k\}_{k \geq 0}\) with the origin being the centre of \(S_0\)
and say that \(S_i\) and \(S_j\) are star adjacent if they share a corner and plus adjacent if they share an edge.
Every square is either vacant or occupied. In this paper, we use the structure of the outermost boundaries derived in
Ganesan (2015) to alternately obtain duality between star and plus connected components in the following sense:
There is a plus connected cycle of vacant squares attached to surrounding the finite star connected component
containing the origin. We also obtain the mutual exclusivity of
left right crossings and top down crossings of star and plus connected components in rectangles.

\vspace{0.1in} \noindent \textbf{Key words:} Star and plus connected components, duality, left right crossings.

\vspace{0.1in} \noindent \textbf{AMS 2000 Subject Classification:} Primary:
60J10, 60K35; Secondary: 60C05, 62E10, 90B15, 91D30.
\end{abstract}

\bigskip

\renewcommand{\theequation}{\thesection.\arabic{equation}}
\setcounter{equation}{0}
\section{Introduction} \label{intro}

Tile \(\mathbb{R}^2\) into disjoint unit squares \(\{S_k\}_{k \geq 0}\) with origin being the centre of~\(S_0.\) Every square in \(\{S_k\}\) is assigned one of the two states, occupied or vacant and the square \(S_0\) containing the origin is always occupied. For \(i \neq j,\) we say that~\(S_i\) and~\(S_j\) are \emph{adjacent} or \emph{star adjacent} if they share a corner between them. We say that \(S_i\) and \(S_j\) are \emph{plus adjacent}, if they share an edge between them. Here we follow the notation of Penrose (2003).

\subsection*{Model Description}
We first discuss star connected components. We say that the square~\(S_i\) is connected to the square~\(S_j\) by a \emph{star connected \(S-\)path} if there is a sequence of distinct squares \((Y_1,Y_2,...,Y_t), Y_l \subset \{S_k\}, 1 \leq l \leq t\) such that~\(Y_l\) is star adjacent to~\(Y_{l+1}\) for all \(1 \leq l \leq t-1\) and \(Y_1 = S_i\) and \(Y_t = S_j.\) If all the squares in \(\{Y_l\}_{1 \leq l \leq t}\) are occupied, we say that~\(S_i\) is connected to~\(S_j\) by an \emph{occupied} star connected \(S-\)path.

Let \(C(0)\) be the collection of all occupied squares in \(\{S_k\}\) each of which is connected to the square~\(S_0\) by an occupied star connected \(S-\)path. We say that \(C(0)\) is the  star connected occupied component containing the origin. Throughout we assume that \(C(0)\) is finite and we study the outermost boundary for finite star connected components containing the origin. By translation, the results hold for arbitrary finite star connected components.



Let \(G_0\) be the graph with vertex set being the set of all corners of the squares of \(\{S_k\}\) in the component~\(C(0)\) and edge set consisting of the edges of the squares of \(\{S_k\}\) in~\(C(0).\) Two vertices in the graph~\(G_0\) are said to be adjacent if they share an edge between them. Two edges in \(G_0\) are said to be adjacent if they share an endvertex between them.

Let \(P = (e_1,e_2,\ldots,e_t)\) be a sequence of distinct edges in \(G_0.\) We say that~\(P\) is a \emph{path} if~\(e_i\) and~\(e_{i+1}\) are adjacent for every \(1 \leq i \leq t-1.\)  Let~\(a\) be the endvertex of~\(e_1\) not common to~\(e_2\) and let~\(b\) be the endvertex of~\(e_t\) not common to~\(e_{t-1}.\) The vertices \(a\) and \(b\) are the \emph{endvertices} of the path~\(P.\)

We say that~\(P\) is a \emph{self avoiding path} if the following three statements hold: The edge \(e_1\) is adjacent only to \(e_2\) and no other~\(e_j, j \neq 2.\) The edge~\(e_t\) is adjacent only to~\(e_{t-1}\) and no other~\(e_j, j \neq t-1.\) For each~\(1 \leq i \leq t-1,\) the edge~\(e_i\) shares one endvertex with~\(e_{i-1}\) and another endvertex with~\(e_{i+1}\) and is not adjacent to any other edge~\(e_j, j\neq i-1,i+1.\)

We say that~\(P\) is a \emph{circuit} if \((e_1,e_2,\ldots,e_{t-1})\) forms a path and the edge~\(e_t\) shares one endvertex with \(e_1\) and another endvertex with~\(e_{t-1}.\) We say that~\(P\) is a \emph{cycle} if \((e_1,e_2,...,e_{t-1})\) is a self avoiding path and the edge~\(e_t\) shares one endvertex with~\(e_1\) and another endvertex with~\(e_{t-1}\) and does not share an endvertex with any other edge~\(e_j, 2 \leq j \leq t-2.\) Any cycle~\(C\) contains at least four edges and divides the plane~\(\mathbb{R}^2\) into two disjoint connected regions. As in Bollob\'as and Riordan~(2006), we denote the bounded region to be the \emph{interior} of~\(C\) and the unbounded region to be the \emph{exterior} of~\(C.\)

We use cycles to define the outermost boundary of star connected components. Let \(e\) be an edge in the graph~\(G_0\) defined above. We say that \(e\) is adjacent to a square \(S_k\) if it is one of the edges of~\(S_k.\) We say that \(e\) is a  \emph{boundary edge} if it is adjacent to a vacant square and is also adjacent to an occupied square of the component~\(C(0).\) Let \(C\) be any cycle of edges in \(G_0.\) We say that the edge~\(e\) is contained in the interior (exterior) of the cycle~\(C\) if at least one endvertex of \(e\) lies in the interior (exterior) of~\(C.\)


We have the following definition.
\begin{Definition} \label{out_def} We say that the edge \(e\) in the graph~\(G_0\) is an \emph{outermost boundary} edge of the component \(C(0)\) if the following holds true for every cycle \(C\) in~\(G_0:\) either \(e\) is an edge in \(C\) or \(e\) belongs to the exterior of \(C.\)

We define the outermost boundary \(\partial _0\) of \(C(0)\) to be the set of all outermost boundary edges of~\(G_0.\)
\end{Definition}

\subsection*{Duality}
To study duality between star and plus connected components, we first define plus connected~\(S-\)cycles and their associated skeletons. Let \(\{J_i\}_{1 \leq i \leq m}\) be a set of squares in the set~\(\{S_k\}.\) We say that the sequence \(L_J = (J_1,...,J_m)\) is a \emph{plus connected \(S-\)path} if for each \(1 \leq i \leq m-1,\) we have that the square~\(J_i\) is plus adjacent i.e., shares an edge, with the square~\(J_{i+1}.\) We say that \(L^+\) is a \emph{plus connected \(S-\)cycle} if \((J_1,\ldots,J_{m-1})\) is a plus connected \(S-\)path and in addition, the square~\(J_m\) is plus adjacent to both~\(J_{m-1}\) and~\(J_1.\)

To define the skeleton of an \(S-\)cycle, we first define the dual graph. Let~\(G\) be the graph with vertex set as the set of corners of the squares in~\(\{S_k\}.\) The edge set of \(G\) is the set of edges of the squares in \(\{S_k\}.\) Let~\(G_d\) be the graph obtained by shifting the graph~\(G\) by \(\left(\frac{1}{2},\frac{1}{2}\right).\) The graph~\(G_d\) also tiles~\(\mathbb{R}^2\) into disjoint unit squares~\(\{W_k\}\) such that \(W_k = S_k + \left(\frac{1}{2},\frac{1}{2}\right).\) We say that the graph \(G\) is the \emph{dual graph} of \(G_d\) and define \(\{W_k\}\) to be dual squares. All quantities related to the dual graph \(G_d\) will also be prefixed with the term dual.

\begin{figure}[tbp]
\centering
\includegraphics[width=2.5in, trim = 120 320 200 220, clip = true]{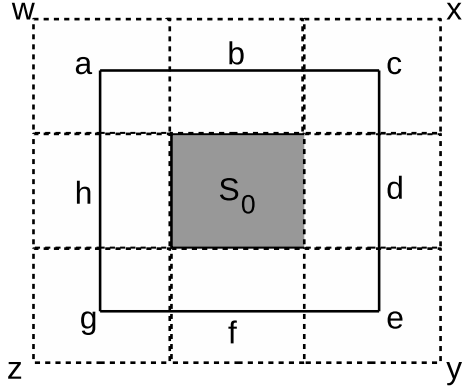}
\caption{Dual skeleton \(SK(L) = abcdefgha\) of the~\(S-\)cycle~\(L\) formed by the dotted squares. }
\label{skel_fig_ill}
\end{figure}

Let \(L_J = (J_1,\ldots,J_m)\) be a plus connected \(S-\)cycle as above. For \(1 \leq k \leq m,\) let \(z_k \in G_d\) be the centre of the square~\(J_k.\) For \(1 \leq i \leq m-1,\) the squares \(J_i\) and \(J_{i+1}\) are plus adjacent; i.e., share an edge between them. Therefore the vertices \(z_{i}\) and \(z_{i+1}\) are adjacent in the dual graph \(G_d;\) i.e., they are connected by a dual edge \(f_i \in G_d.\) Similarly, the vertices~\(z_1\) and~\(z_m\) are also adjacent and are the endvertices of a dual edge~\(f_m \in G_d.\) The resulting sequence of dual edges~\((f_1,f_2,\ldots,f_m)\) is a dual cycle in~\(G_d\) which we denote as~\(SK(L_J).\) We define~\(SK(L_J)\) to be the \emph{dual skeleton} of the plus connected~\(S-\)cycle~\(L_J.\) In Figure~\ref{skel_fig_ill}\((a),\) the sequence of dotted squares with centres~\(a,b,c,d,e,f,g\) and~\(h\) form an \(S-\)cycle~\(L\) of squares. The cycle formed by the vertices~\(abcdefgha\) is a cycle in the dual graph~\(G_d\) and is the dual skeleton~\(SK(L).\) 


Let \(\Lambda_0\) denote the set of all vacant squares in \(\{S_j\}\) that is star adjacent, i.e. shares a corner, with some occupied square in the star connected component~\(C(0).\)
We have the following result.
\begin{Theorem}\label{thm4} Suppose \(C(0)\) is finite. There exists a unique plus connected \(S-\)cycle \(G_{out} = (V_1,...,V_s)\) with the following properties:\\
\((i)\) For every \(i, 1 \leq i \leq s,\) the square~\(V_i \subset \{S_k\}\) is vacant and belongs to~\(\Lambda_0.\)\\
\((ii)\) The outermost boundary of~\(G_{out}\) is a single cycle~\(\partial_G\) consisting only of edges of~\(\{V_i\}_{1 \leq i \leq s}.\) Every vacant square~\(V_i, 1 \leq i \leq s,\) is contained in the interior of~\(\partial_G.\) Every edge in the dual skeleton~\(SK(G_{out})\) is contained in the interior of~\(\partial_G.\)\\
\((iii)\) Every occupied square in the component~\(C(0)\) is contained in the interior of the dual skeleton~\(SK(G_{out}).\) If a vacant square~\(Y \in \Lambda_0\) does not belong to~\(G_{out},\) then~\(Y\) is contained in the interior of~\(SK(G_{out}).\)\\
\((iv)\) If \(F_{out} \neq G_{out}\) is any \(S-\)cycle that satisfies \((i)-(iii)\) above, then the dual skeleton~\(SK(F_{out}) \neq SK(G_{out}).\) Moreover every edge of~\(SK(F_{out})\) either belongs to~\(SK(G_{out})\) or is contained in the interior of~\(SK(G_{out}).\)
\end{Theorem}
The sequence of squares in~\(G_{out}\) form a plus connected \(S-\)cycle of vacant squares surrounding the star connected component~\(C(0).\)

We illustrate the result of Theorem~\ref{thm4} in Figure~\ref{skel_fig_ill} for the simplest case when \(C(0) = S_0,\) the square with origin as centre. The square~\(S_0\) is denoted as the dark grey square and the sequence of surrounding dotted squares form the \(S-\)cycle~\(G_{out}\) of vacant squares. The cycle~\(abcdefgha\) is the dual skeleton~\(SK(G_{out}).\) The dotted square formed by the corners~\(w,x,y\) and~\(z\) forms the outermost boundary~\(\partial_G\) of the \(S-\)cycle~\(G_{out}.\)




For future use we also define star connected \(S-\)cycle. The definition is analogous to plus connected \(S-\)cycles defined above. As before, let \(\{Q_i\}_{1 \leq i \leq n}\) be a set of distinct squares in the set \(\{S_k\}.\) We say that the sequence \(L_Q = (Q_1,...,Q_n)\) is a \emph{star connected \(S-\)path} if for each \(1 \leq i \leq n-1,\) we have that~\(Q_i\) is star adjacent, i.e., shares a corner with~\(Q_{i+1}.\) We say that~\(L_Q\) is a \emph{star connected \(S-\)cycle} if \(n \geq 3,\) the sequence~\((Q_1,\ldots,Q_{n-1})\) forms a star connected \(S-\)path and in addition~\(Q_n\) is also star adjacent to~\(Q_1.\) We consider only star \(S-\)cycles containing at least three squares.


\subsection*{Left right and top bottom crossings in rectangles}
In this subsection, we study the mutual exclusivity of left right and top down crossings in a rectangle. We recall the definition of the graph~\(G\) and the dual graph~\(G_d\) in the previous subsection. The origin is the centre of the square~\(S_0\) and so the dual square \(W_0 = S_0 + \left(\frac{1}{2},\frac{1}{2}\right)\) contains the origin as one of its corners. The rectangle \(R = [0,m] \times [0,n]\) therefore contains exactly~\(mn\) dual squares belonging to~\(\{W_k\}.\) A sequence of distinct squares \(L = (J_1,\ldots,J_m)\) is said to be a \emph{plus connected left right crossing} of \(R\) if the following properties hold.\\
\((a)\) Every square \(J_i \subset \{W_k\}\) is contained in the rectangle~\(R.\)\\
\((b)\) For~\(1 \leq i \leq m-1,\) the square~\(J_i\) is plus adjacent, i.e., shares an edge, with the square~\(J_{i+1}.\)\\
\((c)\) The square~\(J_1\) intersects the left edge of~\(R\) and no other square \(J_i, 2 \leq i \leq m,\) intersects the left edge of \(R.\)\\
\((d)\) The square~\(J_m\) intersects the right edge of~\(R\) and no other square \(J_i, 1 \leq i \leq m-1,\) intersects the right edge of \(R.\)\\
Every square in \(R\) is assigned one of the following two states: occupied or vacant. If every square in a left right crossing~\(L\) is occupied, we say that \(L\) is an \emph{occupied plus connected left right crossing} of the rectangle \(R.\) We denote \(LR^+(R,O)\) and \(LR^+(R,V)\) to be the events that the rectangle~\(R\) contains an occupied and vacant plus connected left right crossing, respectively.

We have a similar definition of plus connected top down crossings and denote \(TD^+(R,O)\) and \(TD^+(R,V)\) to be the events that the rectangle~\(R\) contains an occupied and vacant plus connected top down crossing, respectively. Replacing plus adjacent with star adjacent, we obtain an analogous definition for star connected left right and top down crossings. We have the following result.
\begin{Theorem}\label{thm_lr}We have the following.\\
\((i)\) One of the events \(LR^+(R,O)\) or \(TD^*(R,V)\) always occurs but not both.\\
\((ii)\) One of the events \(LR^*(R,O)\) or \(TD^+(R,V)\) always occurs but not both.
\end{Theorem}
The above result describes the mutual exclusivity of occupied and vacant left right and top down crossings in any rectangle.

The paper is organized as follows: We first prove the main Theorem~\ref{thm4} in Section~\ref{pf4}. Using Theorem~\ref{thm4}, we prove the statement~\((ii)\) of Theorem~\ref{thm_lr} regarding left right crossings in Section~\ref{pf_lr} and the statement~\((i)\) in Section~\ref{pf_lr_i}.







\section{Proof of Theorem~\ref{thm4}}\label{pf4}


\emph{Proof of Theorem~\ref{thm4}}: Let \(\partial_0\) denote the outermost boundary of the star connected component~\(C(0)\) in the graph \(G.\) We recall that the vertex set of~\(G\) is the set of corners of the squares in~\(\{S_k\}\) and the edge set of~\(G\) is the set of edges of the squares in~\(\{S_k\}.\) From Theorem~1 of Ganesan (2015), we have that \(\partial_0 = \cup_{1 \leq i \leq n}C_i\) is a connected union of cycles \(\{C_i\}\) with mutually disjoint interiors and with the property that \(C_i\) and \(C_j\) intersect at most at one point for distinct indices~\(i\) and~\(j.\)

\begin{figure}[ht!]
    \centering
    \begin{subfigure}[t]{0.4\textwidth}
    \centering
        \includegraphics[width=\textwidth, trim = 100 350 200 220, clip = true]{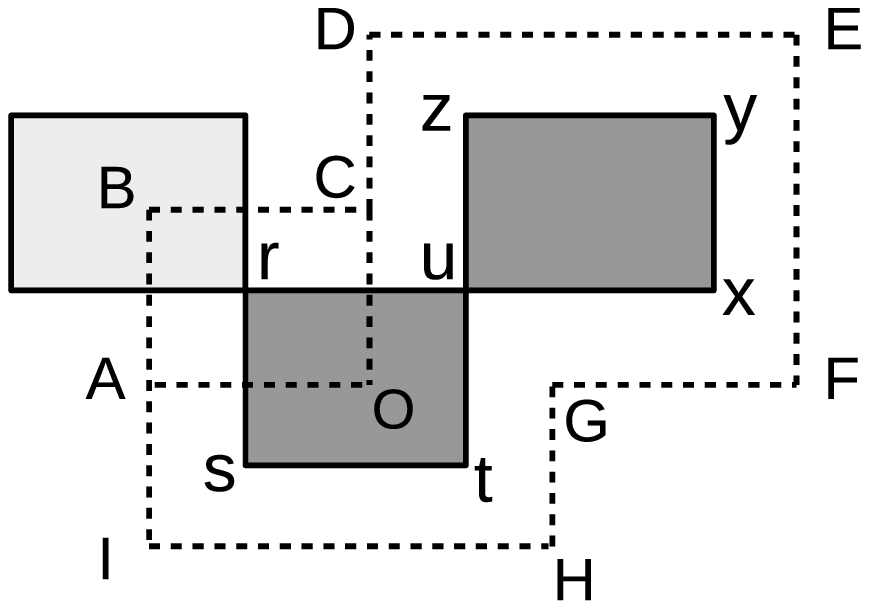}
        \caption{}
     \end{subfigure} 
     ~
     \begin{subfigure}[t]{0.4\textwidth}
     \centering
        \includegraphics[width=\textwidth, trim = 150 270 200 350, clip = true]{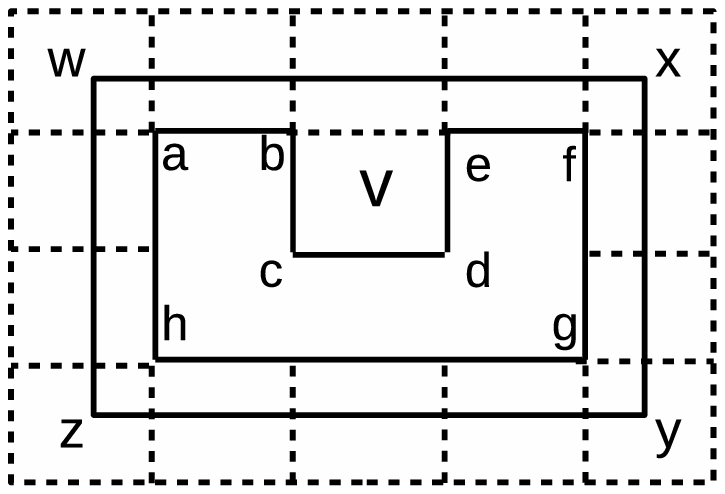}
        \caption{}
    \end{subfigure}

\caption{\((a)\) The component~\(C(0)\) contains two squares denoted by dark grey squares. The dual outermost boundary cycle~\(\partial_V(\partial_0)\)  is denoted by the dotted cycle~\(ABCDEFGHIA.\) Here~\(O\) denotes the origin. \((b)\)~The vacant square marked \(V\) belongs to \(\Lambda\) but is present in the interior of the dual cycle~\(\partial_V(\partial_0)\) denoted by~\(wxyzw.\)}
\label{du_plus_cyc_fig}
\end{figure}

In Figure~\ref{du_plus_cyc_fig}\((a)\), we illustrate the component \(C(0)\) containing two occupied squares denoted by dark grey squares. The outermost boundary~\(\partial_0\) for the component \(C(0)\) is simply the union of the two cycles denoted by~\(urstu\) and~\(uxyzu.\)

To obtain the desired \(S-\)cycle, we use the dual graph \(G_d\) defined prior to the statement of Theorem~\ref{thm4}. We recall that \(G_d\) tiles the plane \(\mathbb{R}^2\) into squares~\(\{W_k\}.\) We define squares in \(\{W_{k}\}\) to be \emph{dual squares} and in general, refer to quantities related to the graph~\(G_d\) with the prefix dual.

We obtain the \(S-\)cycle \(G_{out}\) in three steps. We first obtain a cycle \(C_{dual}\) in the dual graph \(G_d\) containing all the squares of~\(C(0)\) in the interior. In the next step, we show that the dual cycle \(C_{dual}\) is the (dual) skeleton of an \(S-\)cycle \(G_{cyc}\) satisfying properties \((i)-(iii)\) in the statement of the Theorem. In the final step, we merge all the dual skeletons whose \(S-\)cycle satisfies properties~\((i)-(iii),\) to obtain a unique outermost dual skeleton whose corresponding \(S-\) cycle is the desired~\(G_{out}.\)

\subsection*{Obtaining the dual cycle \(C_{dual}\)}
Every vertex \(v\) of the outermost boundary~\(\partial_0\) is the centre of some square \(W_{k(v)} \subset \{W_k\}\) in the dual graph~\(G_d.\) In Figure~\ref{du_plus_cyc_fig}\((a)\), for example, the dual square \(W_{k(r)}\) with centre~\(r\) is denoted by the dotted square~\(AOCB.\)

The union
\begin{equation}\label{cv_par_def}
C_V(\partial_0) = \bigcup_{v \in \partial_0} W_{k(v)}
\end{equation}
is a dual plus connected component since the outermost boundary~\(\partial_0\) is connected in the graph~\(G.\)  Using Theorem~\(2\) of Ganesan~(2015), we therefore have that the dual outermost boundary~\(\partial_V(\partial_0)\) of the dual plus connected component~\(C_V(\partial_0)\) is a single cycle in the dual graph~\(G_d\) containing all (dual) squares of~\(C_V(\partial_0)\) in its interior. Here the dual outermost boundary~\(\partial_V(\partial_0)\) is obtained as follows. Every dual square in \(\{W_k\}\) belonging to~\(C_V(\partial_0)\) is labelled~\(1\) and every dual square in~\(\{W_k\}\) sharing an edge with a dual square in~\(C_V(\partial_0)\) and not belonging to~\(C_V(\partial_0),\) is labelled~\(0.\) We then apply Theorem~\(2\) of Ganesan~(2015) with label \(1\) dual squares as being occupied and label \(0\) dual squares as being vacant.



In Figure~\ref{du_plus_cyc_fig}\((a),\) the dual outermost boundary cycle \(\partial_V(\partial_0)\) is denoted by the dotted cycle~\(ABCDEFGHIA.\)
The dual cycle~\(\partial_V(\partial_0)\) is the desired dual cycle~\(C_{dual}\) and satisfies the following properties. We also refer to the corresponding illustration in Figure~\ref{du_plus_cyc_fig}\((a).\)\\
\((x1)\) Every vertex in the dual cycle~\(\partial_V(\partial_0)\) lies in the exterior of all the cycles belonging to the outermost boundary~\(\partial_0.\) Every edge of~\(\partial_0\) along with its endvertices is contained in the interior of~\(\partial_V(\partial_0).\)\\\\
In Figure~\ref{du_plus_cyc_fig}\((a),\) the vertices \(A,B,\ldots,I\) forming the dual cycle \(\partial_V(\partial_0)\) belong to the exterior of the outermost boundary~\(\partial_0\) formed by the union of the cycles \(urstu\) and \(uxyzu.\)\\\\
\((x2)\) Every occupied square belonging to the component~\(C(0)\) is contained in the interior of the dual cycle~\(\partial_V(\partial_0).\)\\\\
In Figure~\ref{du_plus_cyc_fig}\((a),\) the two dark grey squares forming the component \(C(0)\) lie in the interior of \(\partial_V(\partial_0)\) formed by the dotted cycle \(ABCDEFGHIA.\)\\\\
We recall that \(\Lambda_0\) denotes the set of all vacant squares in~\(\{S_k\}\) sharing a corner with some occupied square in the component~\(C(0).\) Let \(\Lambda \subset \Lambda_0\) be the set of all vacant squares lying in the exterior of the outermost boundary~\(\partial_0\) (i.e., lying in the exterior of all cycles of~\(\partial_0\)) and sharing a vertex with~\(\partial_0.\) \\
\((x3)\) Every vertex in the dual cycle \(\partial_V(\partial_0)\) is the centre of a vacant square in~\(\Lambda.\)\\\\
In Figure~\ref{du_plus_cyc_fig}\((a),\) the light grey square with centre \(B\) belongs to~\(\Lambda.\) Every square with centre in~\(\{A,B,\ldots,I\}\) belongs to~\(\Lambda.\)\\\\
Suppose \(z_1,z_2,\ldots,z_t\) are the vertices of the dual cycle \(\partial_V(\partial_0)\) encountered in that order; i.e., the vertex~\(z_1\) is adjacent to~\(z_2,\) the vertex~\(z_2\) is adjacent to~\(z_3\) and so on. Using property \((x3),\) let \(Y_i\) be the vacant square in \(\Lambda\) with centre~\(z_i,\) for \(1 \leq i \leq t.\)\\
\((x4)\) Every vacant square in \(\Lambda_0\) either belongs to \(\{Y_i\}_{1 \leq i \leq t}\) or lies in the interior of the dual cycle~\(\partial_V(\partial_0).\)\\\\
We remark that the set of squares \(\{Y_i\}_{1 \leq i \leq t}\) need not contain all the vacant squares in \(\Lambda.\) In other words, there may be vacant squares that lie in the exterior of all the cycles of the outermost boundary \(\partial_0\) but lie in the interior of the dual cycle~\(\partial_V(\partial_0).\) This is the case in Figure~\ref{du_plus_cyc_fig}\((b),\) where the outermost boundary \(\partial_0\) is given by the cycle~\(abcdefgh.\) The vacant squares in \(\Lambda\) are shown by the dotted squares and the dual outermost boundary cycle~\(\partial_V(\partial_0)\) is given by the cycle~\(wxyzw.\) The vacant square marked~\(v\) belongs to~\(\Lambda\) but is contained in the interior of~\(\partial_V(\partial_0).\)\\

\emph{Proof of \((x1)-(x4)\)}: We prove the second statement of property \((x1)\) and this also obtains the first statement. Fix an edge~\(e\) in the outermost boundary~\(\partial_0\) with endvertices~\(a_1\) and~\(a_2.\) For illustration, we refer to the edge~\(ru\) in Figure~\ref{du_plus_cyc_fig}\((a)\) with \(a_1 = r\) and \(a_2 = u.\) We recall from (\ref{cv_par_def}) that~\(W_{k(r)}\) and~\(W_{k(u)}\) are the dual squares in~\(C_V(\partial_0)\) with centres \(r\) and \(u,\) respectively. The edge~\(e\) along with its endvertices lie in the interior of~\(W_{k(u)}\cup W_{k(r)}\subset C_V(\partial_0).\) From Theorem~\(2\) of Ganesan~(2015), we have the following property.
\begin{eqnarray}\label{cv_int}
&&\text{The dual outermost boundary cycle~\(\partial_V(\partial_0)\) contains} \nonumber\\
&&\;\;\;\;\;\text{all the (dual) squares of~\(C_V(\partial_0)\) in its interior. }
\end{eqnarray}
Therefore the edge \(e \in \partial_0\) along with its endvertices lies in the interior of the dual cycle~\(\partial_V(\partial_0).\) This proves~\((x1).\)


To prove \((x2),\) we use property \((x1)\) and obtain that the dual cycle~\(\partial_V(\partial_0)\) contains every edge of the outermost boundary~\(\partial_0\) in its interior. Every cycle in~\(\partial_0\) is therefore contained in the interior of \(\partial_V(\partial_0).\) From Theorem~\(1\) of Ganesan~(2015), every occupied square of~\(C(0)\) is contained in the interior of some cycle in~\(\partial_0.\) This proves~\((x2).\)

\begin{figure}[ht!]
    \centering
    \begin{subfigure}[t]{0.4\textwidth}
    \centering
        \includegraphics[width=\textwidth, trim = 100 250 20 20, clip = true]{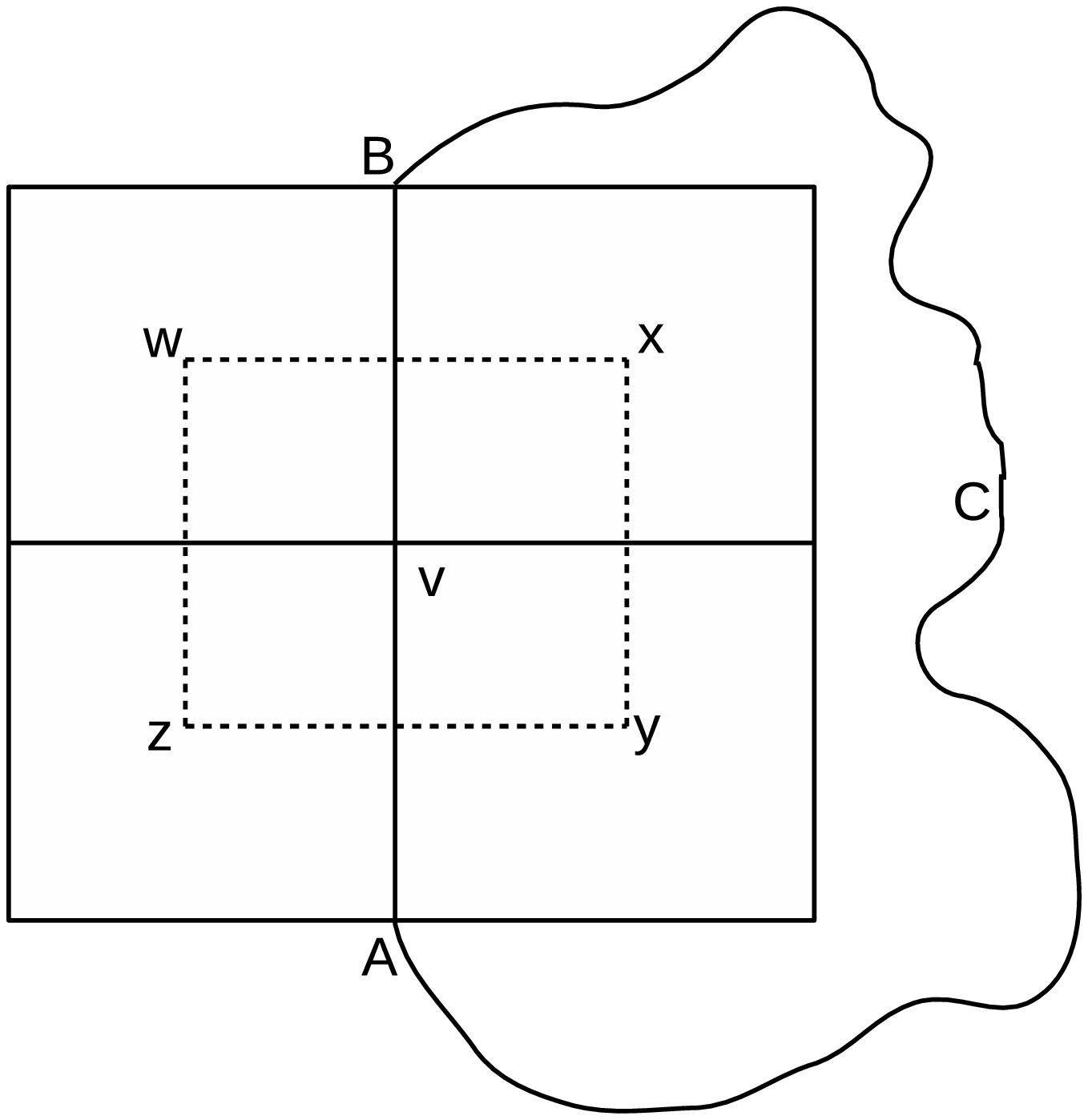}
        \caption{}
     \end{subfigure} 
     ~
     \begin{subfigure}[t]{0.4\textwidth}
     \centering
        \includegraphics[width=\textwidth, trim = 50 350 200 150, clip = true]{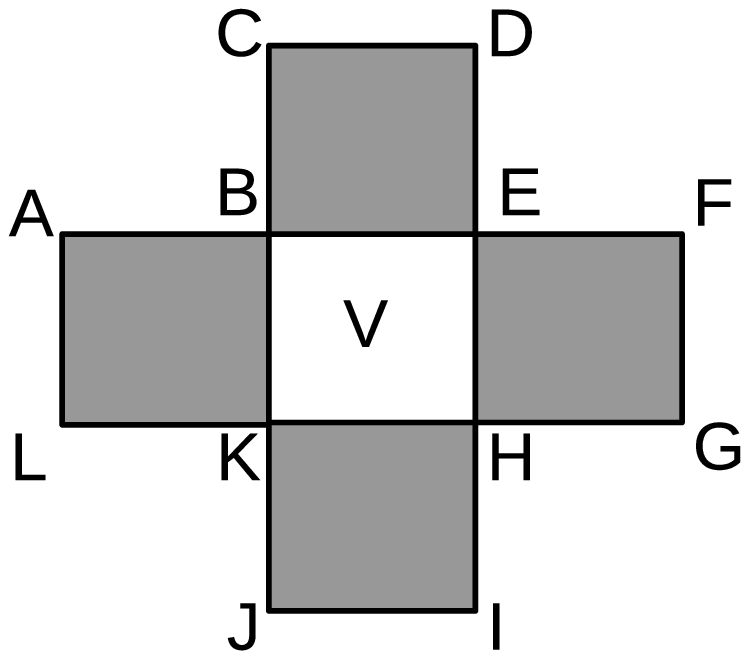}
        \caption{}
    \end{subfigure}

\caption{\((a)\) The dual square \(W_{k(v)}\) with centre \(v\) is depicted by the dotted cycle~\(wxyzw.\) The outermost boundary \(\partial_0\) is denoted by the wavy curve~\(AvBCA.\) \((b)\)~The square \(V\) shares vertices \(B,E,K\) and \(H\) with the outermost boundary \(\partial_0 = ABCEDFGHIJKLA.\) The square \(V\) also lies in the interior of \(\partial_0\) but is vacant.}
\label{j1_fig}
\end{figure}

To prove \((x3)-(x4)\) we use the following property.\\
\((x5)\) Suppose that~\(w\) is a corner of some dual square \(W_{k(v)} \in C_V(\partial_0), v \in \partial_0\) (see (\ref{cv_par_def})). If~\(w\) lies in the exterior of all cycles of~\(\partial_0,\) then~\(w\) is the centre of a vacant square in~\(\Lambda.\)\\
\emph{Proof of \((x5)\)}: Let \(J_1, J_2, J_3\) and \(J_4\) be all the squares in \(\{S_k\}\) containing the vertex~\(v\) as a corner. By construction of the dual graph \(G_d,\) one of these four squares, say \(J_1,\) contains \(w\) as a centre. This is illustrated in Figure~\ref{j1_fig}\((a),\) where the dual square \(W_{k(v)}\) is denoted by the dotted cycle \(wxyzw.\) The square with~\(w\) as centre is \(J_1.\) The outermost boundary~\(\partial_0\) is denoted by the wavy curve~\(AvBCA.\) Since the vertex \(w\) lies in the exterior of all cycles of the outermost boundary~\(\partial_0,\) we have the following property.
\begin{equation}\label{j1_prop}
\text{The square~\(J_1\) lies in the exterior of all cycles of \(\partial_0.\) }
\end{equation}



It remains to see that the square~\(J_1\) is vacant and that at least one of the squares~\(J_2,J_3\) or~\(J_4\) is occupied and belongs to the star connected component~\(C(0).\) This would imply that \(J_1 \in \Lambda.\) We proceed as follows. The vertex~\(v \in \partial_0\) and so at least one of the four edges in the graph~\(G\) containing~\(v\) as an endvertex belongs to the outermost boundary~\(\partial_0.\) We label the edge satisfying the statement above as~\(e_v.\) In Figure~\ref{j1_fig}\((a),\) both the segments \(Av\) and \(vB\) belong to \(\partial_0.\) We fix \(e_v = Av.\)

From Theorem~\(1\) of Ganesan (2015) regarding the outermost boundary~\(\partial_0,\) there are two squares \(J_{a_1},J_{a_2} \in \{J_i\}_{1 \leq i \leq 4}\) containing~\(e_v\) as an edge and exactly one of them, say~\(J_{a_1},\) is occupied and lies in the interior of some cycle of~\(\partial_0.\) The other square~\(J_{a_2}\) is vacant and lies in the exterior of all cycles of~\(\partial_0.\) In Figure~\ref{j1_fig}\((a),\) the square~\(J_{a_1}\) is the square with centre \(y\) and \(J_{a_2}\) is the square with centre~\(z.\)

The square \(J_1\) lies in the exterior of~\(\partial_0\) (see (\ref{j1_prop})) and so we have that \(a_1 \neq 1.\) The square~\(J_{a_1}\) belongs to the  star connected component~\(C(0)\) and shares the vertex~\(v\) with the square~\(J_1.\) If~\(J_{1}\) were also occupied, then \(J_1 \in C(0)\) and would lie in the interior of some cycle of the outermost boundary~\(\partial_0.\) This would contradict (\ref{j1_prop}). Thus \(J_1\) is vacant and this proves \((x5).\)

In the above argument, we need to use the edge \(e_v\) to obtain the occupied and vacant squares, \(J_{a_1}\)  and \(J_{a_2},\) containing \(e_v\) as an edge. This is because, in general, if a square \(A \in \{S_k\}\) lies in the interior of some cycle of \(\partial_0\) and shares a \emph{vertex} with~\(\partial_0,\) it is not necessary that \(A\) is occupied. This is the case in Figure~\ref{j1_fig}\((b),\) where the four dark grey squares form the component~\(C(0).\) The outermost boundary \(\partial_0\) is given by the cycle \(ABCEDFGHIJKLA.\) The vacant square marked \(V\) is contained in the interior of \(\partial_0\) and shares vertices \(B,E,H\) and \(K\) with \(\partial_0.\)\(\qed\)


The property \((x3)\) follows from properties~\((x1)\) and~\((x5)\) since vertices in the dual cycle~\(\partial_V(\partial_0)\) are corners of dual squares in~\(\{W_k\}.\) To see \((x4)\) is true, let \(Z\) be any vacant square in~\(\Lambda_0.\) It suffices to see that the centre~\(q\) of~\(Z\) either belongs to the dual cycle~\(\partial_V(\partial_0)\) or lies in the interior of~\(\partial_V(\partial_0).\) Suppose not and~\(q\) lies in the exterior of~\(\partial_V(\partial_0).\) The vacant square~\(Z\) then lies in the exterior of~\(\partial_V(\partial_0).\) The square~\(Z\) shares a vertex with some occupied square \(J_Z \in C(0).\) Thus some point in \(J_Z\) either belongs to or lies in the exterior of~\(\partial_V(\partial_0).\) Since the occupied square \(J_Z \in C(0)\) is contained in some cycle of the outermost boundary~\(\partial_0,\) we have that some point in the outermost boundary~\(\partial_0\) either belongs to or lies in the exterior of~\(\partial_V(\partial_0).\) This a contradiction to property~\((x1)\) and we obtain~\((x4).\)~\(\qed\)


Let \(L = (Y_1,\ldots,Y_t)\) be the sequence of vacant squares obtained in property \((x3)\) of the dual cycle~\(\partial_V(\partial_0)\) above. We see below that the sequence \(L\) satisfies properties~\((i)-(iii)\) in the statement of the Theorem.

\subsection*{The sequence \(L\) satisfies properties \((i)-(iii)\)}
We use properties \((x2)-(x4)\) above to see that the sequence of squares \(L = (Y_1,\ldots,Y_t)\) is an \(S-\)cycle satisfying the properties~\((i)-(iii)\) in the statement of the theorem. The dual outermost boundary cycle \(\partial_V(\partial_0)\) is connected and since the squares in~\(L\) have their centres in \(\partial_V(\partial_0),\) the sequence~\(L\) is a plus connected \(S-\)cycle. Also we have the following property.
\begin{equation}\label{sk_l}
\text{The dual cycle~\(\partial_V(\partial_0)\) is the dual skeleton~\(SK(L)\) of the \(S-\)cycle~\(L.\)}
\end{equation}

Property~\((i)\) is satisfied by property \((x3)\) above. To see property~\((ii)\) is satisfied we argue as follows. Since the sequence \(L\) is plus connected, we have by Theorem~\(2\) of Ganesan~(2015) that the outermost boundary~\(\partial_0(L)\) of the component~\(L\) satisfies the following properties.\\
\((a1)\) The graph~\(\partial_0(L)\) is a single cycle containing all the squares of~\(L\) in its interior.\\
\((a2)\) Every edge in the outermost boundary~\(\partial_0(L)\) is the edge of some vacant square~\(Y_j \in L.\)\\
The first and the second statements of property~\((ii)\) follow from \((a1)-(a2).\) For the third statement of~\((ii),\) let~\(e\) be an edge of the dual skeleton~\(SK(L)\) of the \(S-\)cycle~\(L.\) The dual edge~\(e\) joins the centres of two squares in~\(L\) and therefore lies in the interior of the union of the squares forming the component~\(L.\) From \((a1),\) we have that the dual edge \(e\) along  with its endvertices lies in the interior of~\(\partial_0(L).\) This proves the second statement of~\((ii).\) 

The second statement of property~\((iii)\) follows from property~\((x4)\) above. The first statement of property~\((iii)\) is true as follows. By property~\((x2)\) of the dual cycle~\(\partial_V(\partial_0),\) every occupied square of the component \(C(0)\) is contained in the interior of~\(\partial_V(\partial_0).\) From property~(\ref{sk_l}) above, we have that~\(\partial_V(\partial_0)\) is the dual skeleton~\(SK(L)\) of the sequence~\(L.\) This proves the first statement of~\((iii).\)



\subsection*{Obtaining the \(S-\)cycle satisfying property \((iv)\)}
\begin{figure}[tbp]
\centering
\includegraphics[width=3.5in, trim= 220 250 90 370, clip=true]{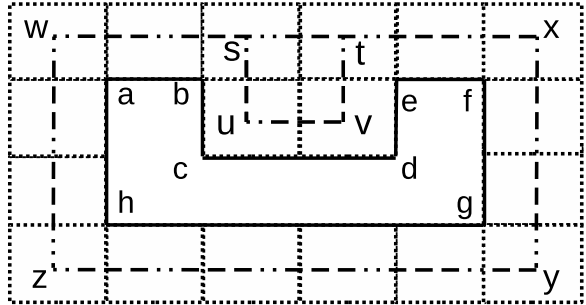}
\caption{Two skeletons \(wsuvtxyzw\) and \(wstxyzw\) for the outermost boundary \(\partial_0 = abcdefgha.\)}
\label{two_skel_fig}
\end{figure}

The \(S-\)cycle \(L\) and the corresponding dual skeleton \(\partial_V(\partial_0)\) satisfy the properties \((i)-(iii)\) in the statement of the Theorem. It is however possible there is more than one dual skeleton satisfying properties \((i)-(iii)\) above. We illustrate this in Figure~\ref{two_skel_fig} where there are two dual skeletons represented by the dotted cycles~\(wstxyzw\) and~\(wsuvtxyzw\) for the outermost boundary~\(\partial_0\) represented by the solid cycle~\(abcdefgha.\) Every dotted square is a vacant square belonging to~\(\Lambda_0;\) i.e., sharing a corner with some occupied square of the component~\(C(0).\)



To obtain a unique \(S-\)cycle that also satisfies property~\((iv),\) we proceed as follows. We perform an iterative procedure starting with the dual cycle~\(\partial_V(\partial_0)\) and obtain an ``outermost" dual cycle~\(K_{out}\) containing~\(\partial_V(\partial_0)\) in its interior. We then argue that~\(K_{out}\) is the dual skeleton of the desired \(S-\)cycle~\(G_{out}\) mentioned in the statement of the Theorem. We recall that a dual skeleton is a cycle in the dual graph~\(G_d\) defined prior to the statement of the Theorem~\ref{thm4}. Let \(K_0 := \partial_V(\partial_0).\) Suppose there is an edge~\(e \in G_d\) with endvertices~\(a,b \in K_0,\) lying in the exterior of~\(K_0.\) Let~\(K_0 = P_0 \cup Q_0\) where~\(P_0\) and~\(Q_0\) are two subpaths of the dual cycle~\(K_0\) with endvertices~\(a\) and~\(b.\)

The unions~\(P_0 \cup \{e\}\) and~\(Q_0 \cup \{e\}\) are both dual cycles and one of the dual cycles, say~\(P_0 \cup \{e\}\) contains~\(K_0\) in its interior. The dual cycle~\(K_1 := P_0 \cup \{e\} \neq K_0\) since the edge~\(e\) lies in the exterior of~\(K_0.\) We set~\(K_1\) to be the dual cycle obtained at the end of the first iteration.

In Figure~\ref{two_skel_fig}, the dual cycle~\(\partial_V(\partial_0) = wsuvtxyzw\) and the edge~\(e\) lying in the exterior of~\(\partial_V(\partial_0)\) is the edge~\(st.\) The dual path~\(P_0 = swzyxt\) and the dual path~\(Q_0 = suvt.\) The dual cycle~\(K_1 = P_0 \cup \{st\} = wstxyzw\) contains the dual cycle~\(\partial_V(\partial_0)\) in its interior.

We repeat the above procedure until there are no exterior edges and denote the final resulting dual cycle as~\(K_{fin}.\)  We have the following properties.\\
\((k1)\) The dual cycle~\(K_{fin} \in {\cal S};\) i.e., the dual cycle~\(K_{fin}\) is the dual skeleton of some \(S-\)cycle~\(E_{fin}\) satisfying
properties~\((i)-(iii)\) in the statement of the Theorem.\\
Let~\({\cal S}\) be the set of all dual skeletons of \(S-\)cycles satisfying the properties~\((i)-(iii)\) in the statement of the Theorem. The set \({\cal S} \neq \emptyset \) since the dual cycle~\(\partial_V(\partial_0) \in {\cal S}.\)\\
\((k2)\) If \(C_{sk} \neq K_{fin}\) is a dual skeleton in~\({\cal S},\) then every dual edge of~\(C_{sk}\) either belongs to~\(K_{fin}\) or is contained in the interior of~\(K_{fin}.\)

Thus the dual skeleton~\(K_{fin}\) also satisfies property~\((iv)\) in the statement of the Theorem.

\emph{Proof of \((k1)-(k2)\)}: The proof of \((k1)\) is true since at the end of every iteration we obtain a dual cycle with the same vertex set as~\(\partial_V(\partial_0)\) and containing~\(\partial_V(\partial_0)\) in its interior. The proof of~\((k2)\) follows from the following property.\\
\((k3)\) If \(C_{sk}\) is any dual skeleton in~\({\cal S},\) then the vertex set of~\(C_{sk}\) is equal to the vertex set of the dual cycle~\(\partial_V(\partial_0).\)

\emph{Proof of \((k3)\)}: We argue by contradiction. Let \(L_{sk}\) be the \(S-\)cycle whose dual skeleton is~\(C_{sk}.\) Suppose there is a vertex~\(v \in C_{sk}\) lying in the exterior of~\(\partial_V(\partial_0).\) The square~\(Y_v \in \{S_k\}\) with centre~\(v\) is vacant and belongs to~\(\Lambda_0\) since~\(L_{sk}\) satisfies properties~\((i)-(iii).\) We recall that~\(\Lambda_0\) is the set of all vacant squares in~\(\{S_k\}\) sharing a corner with some occupied square of~\(C(0).\) The square~\(Y_{v} \in \Lambda_0\) lies in the exterior of the dual cycle~\(\partial_V(\partial_0),\) contradicting the fact that~\(\partial_V(\partial_0)\) satisfies property~\((ii).\) An analogous argument obtains that there cannot be a vertex of~\(\partial_V(\partial_0)\) lying in the exterior of~\(C_{sk}.\) \(\qed\)

Using property~\((k3),\) we obtain~\((k2)\) as follows. If \(C_{sk} \neq K_{fin},\) then there is an edge~\(e_{sk}\) of~\(C_{sk}\) not present in~\(K_{fin}.\) From property~\((k3),\) both the endvertices of~\(e_{sk}\) belong to~\(\partial_V(\partial_0)\) and therefore belong to~\(K_{fin}.\) If~\(e_{sk}\) lies in the exterior of~\(K_{fin},\) the above iterative procedure would not have terminated. This proves~\((k2).\)\(\qed\)

\section{Proof of~\((ii)\) in Theorem~\ref{thm_lr}}\label{pf_lr}



We first argue the following.
\begin{equation}\label{lr_td_int}
\text{The events \(LR^*(R,O)\) and \(TD^+(R,V)\) cannot occur simultaneously.}
\end{equation}
\emph{Proof of~(\ref{lr_td_int})}: Suppose not and suppose \(\Pi^*_R\) is an occupied star connected left right crossing and suppose \(\Pi_T\) is a vacant plus connected top down crossing.

We first construct two paths \(\Gamma_T\) and \(\Gamma^*_R\) using the crossings \(\Pi_T\) and \(\Pi^*_R,\) respectively. We use the paths \(\Gamma_T\) and \(\Gamma^*_R\) to arrive at a contradiction.

\emph{\underline{Construction of \(\Gamma_T\)}}: The crossing \(\Pi_T = (A_1,\ldots,A_t)\) is a plus connected component with the square~\(A_1 \in \{W_k\}\) intersecting the top edge of the rectangle~\(R\) and the square~\(A_t \in \{W_k\}\) intersecting the bottom edge of the rectangle \(R.\) We recall the definition of the dual graph \(G_d\) and the associated dual squares \(\{W_k\}\) prior to the statement of Theorem~\ref{thm4}. For the purposes of this proof, we call the graph~\(G\) as the dual graph to the graph~\(G_d.\)

We now extend the crossing \(\Pi_T\) slightly as follows. Let \(A_0\) be the square lying above~\(A_1\) and sharing an edge with~\(A_1.\) Similarly, let \(A_{t+1}\) be the square lying below~\(A_t\) and sharing an edge with~\(A_t.\) The square \(A_0\) lies in the exterior of the rectangle \(R\) and the bottom edge of \(A_0\) is contained in the top edge of \(R.\) Similarly, the square \(A_{t+1}\) lies below the rectangle~\(R\) and intersects the bottom edge of \(R.\) For \(0 \leq i \leq t+1,\) let~\(z_i\) be the centre of the square~\(A_i.\) The path \(\Gamma_T = (g_0,g_2,\ldots,g_{t})\) formed by the vertices \(z_0,\ldots,z_{t+1}\) is a path in the dual graph~\(G.\) For \(1 \leq j \leq t,\) the dual edge~\(g_{j}\) has endvertices \(z_{j}\) and~\(z_{j+1}.\)

\emph{\underline{Construction of \(\Gamma^*_R\)}}: Let \(\partial_R\) be the outermost boundary of the star connected crossing~\(\Pi^*_R.\) We have the following property.\\
\((q)\) There is a path \(\Gamma^*_R  = (f_1,\ldots,f_w)\subset \partial_R\) such that the edge~\(f_1\) intersects the left edge of the rectangle~\(R,\) the edge~\(f_w\) intersects the right edge of \(R\) and every edge \(f_i, 1 \leq i \leq w\) is contained in \(R.\)\\
We say that an edge~\(e\) is contained in the rectangle~\(R\) if no endvertex of~\(e\) lies in the exterior of~\(R.\)

\emph{Proof of \((q)\)}: The crossing \(\Pi^*_R\) is a star connected component and from Theorem~\(1\) of Ganesan~(2015), the outermost boundary~\(\partial_R\) of the component~\(\Pi^*_R\) is a connected union of cycles \(\cup_{i=1}^{w}C_i\) with mutually disjoint interiors. The edges of~\(\partial_R\) are edges of squares in the crossing~\(\Pi^*_R\) and so every edge in~\(\partial_R\) is contained in the rectangle~\(R.\) Since~\(\Pi^*_R\) is also a left right crossing, it contains two squares~\(Q_1\) and~\(Q_2\) such that the square \(Q_1\) intersects the left edge and the square \(Q_2\) intersects the right edge of the rectangle~\(R.\)

The square~\(Q_1\) is the unique square in~\(\Pi^*_R\) that intersects the left edge of \(R.\) Also since \(Q_1\) is contained in the interior of some cycle \(C_i \in \partial_R,\) the bottom edge \(e_b \in Q_1\) belongs to~\(C_i.\) The edge \(e_b\) intersects the left edge of~\(R.\) Similarly, the bottom edge \(f_b\) of the square \(Q_2\) intersects the right edge of \(R\) and also belongs to~\(\partial_R.\) Since~\(\partial_R\) is connected, the edges~\(e_b\) and~\(f_b\) are connected by a path of edges in~\(\partial_R\) and we obtain~\((q).\)\(\qed\)

We obtain the desired contradiction using the paths~\(\Gamma_T\) and~\(\Gamma^*_R.\) By construction, the paths~\(\Gamma_T\) and~\(\Gamma^*_R\) intersect in the following sense. There are edges~\(f_k \in \Gamma^*_R\) and~\(g_j \in \Gamma_T\) such that~\(f_k\) and~\(g_j\) intersect. We consider two possible cases depending on whether \((a)\) the index \(j \in \{0,t\}\) or \((b)\) \(1 \leq j \leq t-1.\) For case \((b),\) the argument is as follows. The dual edge~\(f_k\) is the edge common to the dual squares~\(A_j\) and~\(A_{j+1}\) with centres~\(z_j\) and~\(z_{j+1},\) respectively. Since \(1 \leq j \leq t-1,\) both the squares~\(A_j\) and~\(A_{j+1}\) belong to~\(\Pi_T\) and so are both vacant. But the edge~\(f_k\) belongs to the outermost boundary~\(\partial_R\) and so we have from Theorem~\(1\) of Ganesan~(2015), that one of the squares in~\(\{A_j,A_{j+1}\}\) belongs to~\(\Pi^*_R\) and is therefore occupied, a contradiction. This proves~(\ref{lr_td_int}) for case \((b).\)

For~case~\((a),\) we prove for the case when \(j =0.\) An analgous argument holds for the case \(j = t.\) If \(j =0,\) the the edge \(f_k\) is necessarily contained in the top edge of \(R\) and the topmost square \(A_1 \in \Pi_T\) contains \(f_k\) as the top edge. The square \(A_1\) is vacant since \(A_1 \in \Pi_T.\) But by Theorem~\(1\) of Ganesan (2015), one of the squares containing \(f_k\) as an edge belongs to \(\Pi^*_R\) and is necessarily occupied. The square \(A_0\) lying above~\(A_1\) lies in the exterior of the rectangle~\(R\) and therefore does not belong to~\(\Pi^*_R.\) This means that~\(A_1\) must be occupied, a contradiction.~\(\qed\)

\begin{figure}[tbp]
\centering
\includegraphics[width=3.5in, trim= 20 320 90 120, clip=true]{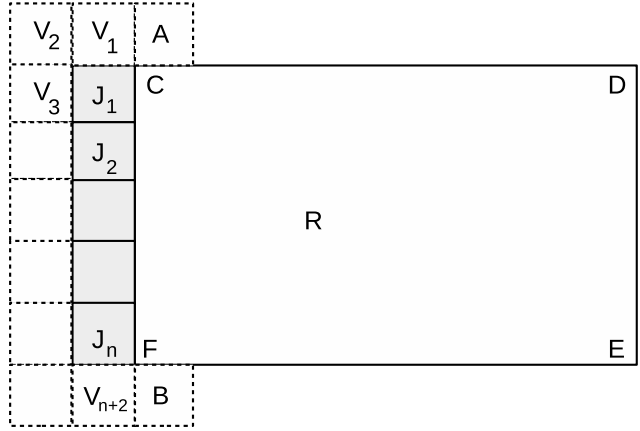}
\caption{The label \(1\) squares \(\{J_i\}_{1 \leq i \leq n}\) and the surrounding label \(0\) squares \(\{V_j\}_{1 \leq j \leq n+2}.\) The rectangle~\(R\) is denoted by \(CDEF.\)}
\label{wj_fig}
\end{figure}

In the rest of the proof, we provide the procedure to obtain a vacant plus connected top down crossing if there is no occupied star connected left right crossing. Let \(R_{left}\) denote the left edge of the rectangle~\(R = [0,m] \times [0,n].\) Let \(\{J_i\}_{1 \leq i \leq n} \subset \{W_k\}\) be the~\(n\) (dual) squares lying in the exterior of~\(R\) and intersecting the left edge of~\(R.\) Unless specifically mentioned, we suppress the notation dual in the rest of the proof for simplicity. If~\((x_i,y_i)\) is the centre of the square~\(J_i,\) then~\(x_i = x_1\) for all~\(1 \leq i \leq n.\) We assume that \(y_1 > y_2 > \ldots > y_n;\) i.e., the squares are indexed in the decreasing of the \(y-\)coordinates of their centre. This is illustrated in Figure~\ref{wj_fig}, where the rectangle \(R\) is denoted by \(CDEF.\)

Label all squares~\(\{J_i\}_{1 \leq i \leq n}\) to be~\(1.\) We label the other squares in~\(\{W_k\}\) using the following procedure. Let \(W \in \{W_k\}\) be any square. Suppose \(W\) is occupied and contained in the interior of the rectangle~\(R.\) We recall the definition of star connected \(S-\)path prior to the statement of Theorem~\ref{thm4} and suppose that~\(W\) is connected to some square~\(J_i, 1 \leq i \leq n,\) by a star connected path of the form \(P = (W,A_1,A_2,\ldots,A_h,J_i),\) where each \(A_j, 1 \leq j \leq h\) is an occupied square in the interior of the rectangle~\(R.\) We label~\(W\) as~\(1.\) The union of all the label \(1\) squares, \(C_{left},\) is a star connected component containing the squares~\(\{J_i\}_{1 \leq i \leq n}.\)

If the square \(W \in \{W_k\}\) shares a corner with some label~\(1\) square (belonging to~\(C_{left}\)) and is not labelled, we label~\(W\) as~\(0.\) (The square~\(W\) may also be present in the exterior of the rectangle \(R.\)) We have the following properties.
\begin{eqnarray}
&&\text{If \(W \in \{W_k\}\) lies in the interior of the rectangle~\(R\)}\nonumber\\
&&\;\;\;\;\;\text{and has label~\(0,\) then~\(W\) is vacant.}\label{w_int}
\end{eqnarray}
\begin{eqnarray}
&&\text{The rectangle~\(R\) contains an occupied left right star connected crossing }\nonumber\\
&&\;\;\;\;\text{if and only if the component~\(C_{left}\) contains a (label \(1\)) square }\nonumber\\
&&\;\;\;\;\;\text{intersecting the right edge of~\(R.\)}\label{lr_c_left}
\end{eqnarray}




We apply Theorem~\ref{thm4} to the star connected component \(C_{left}\) with occupied replaced by label~\(1\) and vacant replaced by label~\(0.\) Let~\(L\) be the plus connected \(S-\)cycle of label~\(0\) squares surrounding~\(C_{left}\) and satisfying the corresponding properties~\((i)-(iv).\)


We have the following properties. \\
\((n1)\) Let \(V_1,\ldots,V_{n+2}\) be the squares lying to the left of the line \(x = 0\) and sharing edges with squares in~\(\{J_i\}_{1 \leq i \leq n},\) as illustrated in Figure~\ref{wj_fig}. The square \(V_{i}, 1 \leq i \leq n+2,\) has label~\(0\) and belongs to the sequence~\(L.\) The sequence \((V_1,\ldots,V_{n+2})\) is a subsequence of consecutive squares in~\(L.\)\\\\
Let \(L_I = L \setminus (V_1,\ldots,V_{n+2}) = (Y_{1},\ldots,Y_h).\)\\
\((n2)\) The square~\(Y_{1}\) intersects the top edge of~\(R\) and shares an edge with the square~\(V_1.\) Similarly, the square~\(Y_h\) intersects the bottom edge of~\(R\) and shares an edge with~\(V_{n+2}.\)\\
In Figure~\ref{wj_fig}, the square~\(Y_{1}\) is labelled~\(A\) and the square~\(Y_h\) is labelled~\(B.\)\\\\
Using property \((n2),\) let \(j_1\) be the largest index \(j \geq 1\) such that the square \(Y_j \in L_I\) intersects the top edge of~\(R\) and let~\(j_2\) be the smallest index~\(k \leq h\) such that~\(Y_k \in L_I\) intersects the bottom edge of~\(R.\) Define the sequence \(L_{II} = (Y_{j_1},Y_{j_1+1},\ldots,Y_{j_2}).\)\\
\((n3)\) The sequence \(L_{II}\) is a plus connected \(S-\)path (see definition prior to statement of Theorem~\ref{thm4}). If there is no occupied star connected left right crossing of~\(R,\) then every square in the sequence~\(L_{II}\) lies in the interior of the rectangle~\(R\) and is vacant.\\

The sequence~\(L_{II}\) is the desired vacant plus connected top down crossing of~\(R.\)


\emph{Proof of \((n1)-(n3)\)}: The first half of property~\((n1)\) is true by the labelling procedure. We prove the second half of~\((n1)\) as follows. We prove for the square \(V_1\) and extend the proof for \(V_i, 2 \leq i \leq n+2.\) Let~\(\partial_0(L)\) be the outermost boundary of the sequence~\(L.\) From property~\((ii)\) of Theorem~\ref{thm4}, we have that~\(\partial_0(L)\) is a single cycle consisting of edges of the squares in~\(L\) and from property~\((iii),\) we have that the square~\(J_1\) is contained in the interior of the cycle~\(\partial_0(L).\)

Suppose that the square~\(V_1\) is not present in the sequence~\(L.\) The square~\(V_1\) has label~\(0\) and shares an edge with the label~\(1\) square~\(J_1.\) From properties~\((ii)-(iii)\) of Theorem~\ref{thm4}, we have that~\(V_1 \in \Lambda_0\) is contained in the interior of the outermost boundary cycle~\(\partial_0(L).\) Here~\(\Lambda_0\) is the set of all label~\(0\) squares sharing a corner with some label~\(1\) square in~\(C_{left}.\)

\begin{figure}[tbp]
\centering
\includegraphics[width=3.5in, trim= 20 320 90 120, clip=true]{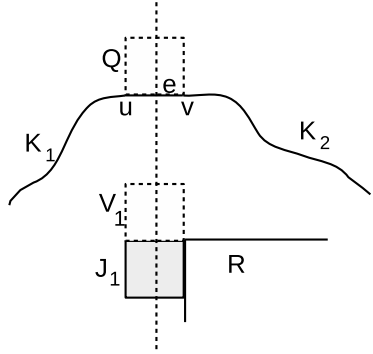}
\caption{The edge \(e\) with endvertices~\(u\) and~\(v\) lies above the label~\(0\) square~\(V_1\) and is also the edge of some label~\(0\) square~\(Q.\) The long dotted line denotes the line~\(x = x_1\) and the wavy lines \(K_1\) and \(K_2\) are parts of the outermost boundary cycle \(\partial_0(L).\)}
\label{e_int}
\end{figure}

We now obtain a contradiction by showing that there is a labelled square above~\(V_1.\) Since the square~\(J_1\) lies in the interior of the outermost boundary cycle~\(\partial_0(L),\) some edge~\(e \in \partial_0(L)\) intersects the line~\(x = x_1\) at~\((x_1,y_e).\) We recall that~\((x_1,y_1)\) and~\((x_1,y_1+1)\) are the centres of the squares~\(J_1\) and~\(V_1,\) respectively. Since the square~\(V_1\) also lies in the interior of~\(\partial_0(L),\) the edge~\(e\) does not belong to the square~\(V_1\) and lies above~\(V_1\) (see Figure~\ref{e_int}).
From property \((ii)\) of Theorem~\ref{thm4}, we have the edge \(e\) is the edge of a label~\(0\) square~\(Q.\) The square~\(Q\) lies above square~\(V_1\) and this is a contradiction since there is no labelled square above the square~\(V_1.\)

From the above discussion, we obtain that~\(V_1 \in L\) and this proves the first statement in~\((n1)\) for~\(V_1.\) An analogous proof holds for the bottom square~\(V_{n+2}\) and the squares \(V_i, 3 \leq i \leq n.\) The squares~\(V_2\) and~\(V_{n+1}\) also belong to the sequence~\(L\) and we give below an argument for~\(V_2\) and an analogous argument holds for~\(V_{n+1}.\) The square~\(V_2\) shares only a corner with the square~\(J_1\) and does not share a corner with any other square in~\(\{J_i\}_{2 \leq i \leq n}.\) The sequence~\(L\) is plus connected and therefore the label~\(0\) square~\(V_1\) is plus adjacent to two label~\(0\) squares. The square sharing the top edge of~\(V_1\) does not have any label. The square sharing the bottom edge of~\(V_1\) is the square~\(J_1\) which has label~\(1.\) Therefore the square marked~\(A\) in Figure~\ref{wj_fig} and the square~\(V_2\) have label~\(0\) and belong to~\(L.\)

To see that~\((V_1,V_2,\ldots,V_{n+2})\) is a subsequence of the sequence~\(L,\) we have by construction that~\(L\) is a plus connected \(S-\)cycle of label~\(0\) squares. For \(2 \leq i \leq n-1,\) the only label~\(0\) squares plus adjacent, i.e., sharing an edge with~\(V_{i}\) are~\(V_{i-1}\) and~\(V_{i+1}.\) This proves~\((n1).\) From the discussion in the previous paragraph, we also have that the square marked~\(A\) is the square~\(Y_{1} \in L\) sharing an edge with the square~\(W_1\) and the square marked~\(B\) is the square~\(Y_h \in L\) sharing an edge with the square~\(W_{n+2}.\) This proves~\((n2).\)


The first statement in property~\((n3)\) is true since~\(L\) is a plus connected \(S-\)cycle and~\(L_{II}\) is a subsequence of consecutive squares in~\(L.\) To see that the second statement of~\((n3)\) is true, suppose that some square~\(Y_j \in L_{II}\) lies in the exterior of the rectangle~\(R.\) If~\(Y_j\) intersects the top edge of~\(R,\) then since~\(Y_{j+1}\) is plus adjacent to~\(Y_j,\) the square~\(Y_{j+1}\) necessarily contains either the top, left or the bottom edge of~\(Y_j;\) i.e., the square~\(Y_{j+1}\) cannot lie above~\(Y_j.\) This is because all labelled squares either lie in the interior of the rectangle~\(R\) or intersect the boundary of~\(R.\) Thus the square~\(Y_{j+1}\) also intersects the top edge of~\(R,\) a contradiction to the definition of the index~\(j_1.\)

An analogous arugment as in the previous paragraph also obtains that the square~\(Y_j\) cannot intersect the bottom edge of~\(R.\) If~\(Y_j\) intersects the left edge of the rectangle~\(R,\) then \(Y_j \in \{J_i\}_{1 \leq i \leq n}.\) This is a contradiction since~\(Y_j\) has label~\(0\) and every square~\(J_i, 1 \leq i \leq n,\) has label~\(1.\) If~\(Y_j\) intersects the right edge of~\(R,\) then there is a label~\(1\) square \(Q_j \in C_{left}\) sharing a corner with~\(Y_j.\) The square \(Q_j\) lies in the interior of the rectangle \(R\) and touches the right edge of~\(R.\) This contradicts property~(\ref{lr_c_left}). This proves \((n3).\)~\(\qed\)


\section{Proof of~\((i)\) in Theorem~\ref{thm_lr}}\label{pf_lr_i}
\emph{Proof of \((i)\)}: As in the proof of \((ii),\) we have that the events \(LR^+(R,O)\) and \(TD^*(R,V)\) cannot occur together.

As in the proof of \((ii),\) let \(J_i ,1 \leq i \leq n,\) be the squares intersecting the left edge of the rectangle~\(R\) and lying in the exterior of~\(R.\) We again label all squares \(\{J_i\}_{1 \leq i \leq n}\) to be \(1.\) We label the other squares using the following procedure.



Label all squares \(\{J_i\}_{1 \leq i \leq n}\) to be \(1.\) We label the other squares in~\(\{W_k\}\) using the following procedure. Let \(W \in \{W_k\}\) be any square. Suppose \(W\) is occupied and contained in the interior of the rectangle~\(R.\) We recall the definition of plus connected \(S-\)path prior to the statement of Theorem~\ref{thm4} and suppose that~\(W\) is connected to some square~\(J_i, 1 \leq i \leq n,\) by a plus connected path of the form \(P = (W,A_1,A_2,\ldots,A_h,J_i),\) where each \(A_j, 1 \leq j \leq h\) is an occupied square in the interior of the rectangle~\(R.\) We label~\(W\) as~\(1.\) The union of all the label \(1\) squares, \(C^+_{left},\) is a plus connected component containing the squares~\(\{J_i\}_{1 \leq i \leq n}.\)

If the square~\(W \in \{W_k\}\) shares an edge with some label~\(1\) square (belonging to~\(C^+_{left}\)) and is not labelled, we label~\(W\) as~\(0.\) (The square~\(W\) may also be present in the exterior of the rectangle~\(R.\)) Analogous to~(\ref{w_int}) and~(\ref{lr_c_left_plus}), we have the following properties.
\begin{eqnarray}
&&\text{If a square~\(W \in \{W_k\}\) lies in the interior of the rectangle~\(R\)}\nonumber\\
&&\;\;\;\;\;\text{and has label~\(0,\) then~\(W\) is vacant.}\label{w_int_plus}
\end{eqnarray}
\begin{eqnarray}
&&\text{The rectangle~\(R\) contains an occupied left right plus connected crossing }\nonumber\\
&&\;\;\;\;\text{if and only if the component~\(C^+_{left}\) contains a (label~\(1\)) square }\nonumber\\
&&\;\;\;\;\;\text{intersecting the right edge of~\(R.\)}\label{lr_c_left_plus}
\end{eqnarray}

Let \(\partial^+_{left}\) be the outermost boundary of the plus connected~\(C^+_{left}.\)  From Theorem~\(2\) of Ganesan~(2015), we have that~\(\partial^+_{left}\) is a single cycle containing all the squares of~\(C^+_{left}\) in its interior.

The proof consists of three steps. In the first step of the proof, we extract a set \(\{A_j\}\) of a label~\(0\) squares sharing edges with~\(\partial^+_{left}.\) In the second step of the proof, we merge the squares~\(\{A_j\}\) with~\(\partial^+_{left}\) iteratively to obtain a final cycle~\(D_{fin}.\) In the final step of the proof, we construct the desired top down crossing of vacant squares by exploring the cycle~\(D_{fin}\) edge by edge.


\subsection*{Extracting the label~\(0\) squares \(\{A_j\}\) attached to \(\partial^+_{left}\)}
\begin{figure}[tbp]
\centering
\includegraphics[width=3.5in, trim= 40 340 90 120, clip=true]{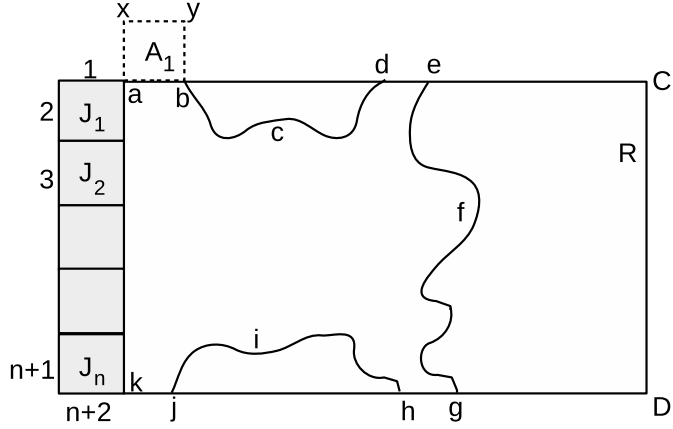}
\caption{The label~\(1\) squares~\(\{J_i\}_{1 \leq i \leq n}\) lying outside the rectangle~\(R\) are shown in light grey. The rectangle~\(R\) is denoted by~\(CDka.\) The outermost boundary cycle~\(\partial^+_{left}\) of the component~\(C^+_{left}\) is the union of the path of edges~\(Q = (1,2,3,\ldots,n+1,n+2)\) and the wavy path~\(abcdefghijk.\) }
\label{wj_fig_plus}
\end{figure}

We first enumerate the relevant properties of the outermost boundary cycle~\(\partial^+_{left}.\) \\
\((q1)\) For \(1 \leq i \leq n+2,\) let \(h_i\) be the edge denoted by the letter \(i\) in Figure~\ref{wj_fig_plus}. The path \(Q = (h_1,\ldots,h_{n+2})\) with endvertices \(a\) and \(k\) is a subpath of the cycle \(\partial^+_{left}.\)\\\\

Let \(Q_I = \partial^+_{left} \setminus Q = (g_{1},\ldots,g_b).\) In Figure~\ref{wj_fig_plus}, the path \(Q_I\) is the wavy path~\(abcdefghijk.\)\\\\
\((q2)\) The path~\(Q_I\) contains~\((0,0)\) and~\((0,n)\) as endvertices. The edge~\(g_{1}\) contains~\((0,n)\) as an endvertex and therefore intersects the top edge of the rectangle~\(R.\) The edge~\(g_b\) contains~\((0,0)\) as an endvertex and intersects the bottom edge of~\(R.\) Every edge of~\(Q_I\) is contained in~\(R\) in the following sense: If~\(u\) is an endvertex of an edge~\(e \in Q_I,\) then~\(u\) either belongs to the boundary of~\(R\) or is contained in the interior of~\(R.\)\\\\
In Figure~\ref{wj_fig_plus}, the vertex~\(k = (0,0),\) the vertex~\(a = (0,n),\) the edge~\(g_1 = ab\) and the edge \(g_b = kj.\) Every edge in the wavy path~\(Q_I = abcdefghijk\) is contained in~\(R.\)\\\\
From Theorem~\(1\) of Ganesan (2015), we have that every edge \(g_j \in Q_I, 1 \leq j \leq b,\) is the edge of a label~\(1\) square~\(K_j\) belonging to the component~\(C^+_{left}.\) The edge~\(g_j\) is also the edge of a label~\(0\) square~\(A_j\) lying in the exterior of~\(\partial^+_{left}.\)\\\\
\((q3)\) The label~\(1\) square \(K_j\) containing the edge \(g_j\) is contained in the interior of the outermost boundary cycle~\(\partial^+_{left}.\) The label~\(0\) square \(A_j\) lies to the right of the line \(x=0\) in the following sense: If \(v \in A_j\) is a corner, then either~\(v\) belongs to the line~\(x = 0\) or lies to the right of the line~\(x = 0.\)\\\\
In~Figure~\ref{wj_fig_plus}, the line \(x= 0\) is the infinite line containing the segment~\(ak.\)


\emph{Proof of \((q1)-(q3)\)}: The proof \((q1)\) is analogous to the proof~\((n1)\) and is obtained by contradiction assuming that~\(h_1 \notin \partial^+_{left}.\) The cycle~\(\partial^+_{left}\) then intersects the line~\(x = x_1\) at some edge lying above~\(h_1.\) We recall that~\((x_1,y_1)\) is the centre of the label~\(1\) square~\(J_1.\) Arguing as in the proof of~\((n1),\) we then obtain~\((q1).\)

The first statement of property~\((q2)\) is true since the path~\((h_1,\ldots,h_{n+2})\) contains~\((0,n)\) and~\((0,0)\) as endvertices (see Figure~\ref{wj_fig_plus}). For the second statement of~\((q2),\) we consider an edge \(g_j \in Q_I.\) From the first statement of~\((q3)\) we have that~\(g_j\) is the edge of a label~\(1\) square~\(K_j \in C^+_{left}.\) If~\(K_j\) lies in the interior of the rectangle~\(R,\) then the second statement of~\((q2)\) is true.
We also have the following property.
\begin{eqnarray}
&&\text{If \(K_j  = J_i\) for some \(1 \leq i \leq n,\) then the edge \(g_j\)}\nonumber\\
&&\;\;\;\;\text{is necessarily the right edge of~\(J_i\)}\label{g_j_prop}
\end{eqnarray}
Therefore the second statement of~\((q2)\) is again true.\\
\emph{Proof of (\ref{g_j_prop})}: For illustration, we suppose \(K_j = J_1\) and an analogous argument holds for all~\(J_i.\) The edge~\(g_j\) cannot be the left edge (marked~\(2\) in Figure~\ref{wj_fig_plus}) since the edge marked~\(2\) belongs to the subpath~\(Q = \partial^+_{left} \setminus Q_I.\) For the same reason, the edge~\(g_j\) cannot be the top edge of~\(J_1\) (marked~\(1\) in Figure~\ref{wj_fig_plus}). If the height of the rectangle \(n = 1,\) then the bottom edge of~\(J_1\) also belongs to~\(Q\) and cannot be equal to~\(g_j.\)

If the height of the rectangle \(n \geq 2,\) then the bottom edge of~\(J_1\) lies in the interior of~\(\partial^+_{left}.\) This is true since the edges marked~\(2\) and~\(3\) both belong  to the subpath~\(Q\) of the cycle~\(\partial^+_{left}\) and so the bottom edge of~\(J_1\) cannot belong to~\(\partial^+_{left}.\) Thus again in this case the bottom edge cannot be equal to~\(g_j.\) This implies that~\(g_j\) is the right edge of~\(J_i.\) This proves~(\ref{g_j_prop}) and the property~\((q2).\)~\(\qed\)

The final statement in property~\((q3)\) is true since by property~\((q2),\) every edge of~\(Q_I\) is contained in the rectangle~\(R.\) Therefore no endvertex of \(g_j\) lies to the left of the line \(x= 0.\) If the edge~\(g_j\) contains a endvertex to the right of the line \(x = 0,\) then the corresponding label~\(0\) square \(A_j\) lies to the right of the line \(x = 0.\) Suppose now that~\(g_j\) is contained in the line \(x=0.\) By Theorem~\(1\) of Ganesan~(2015), the edge \(g_j\) is common to a label~\(1\) square and a label~\(0\) square. We therefore necessarily have that \(K_j = J_i\) for some \(1 \leq i \leq n.\) The label~\(0\) square~\(A_j\) again lies to the right of the line~\(x=0.\)~\(\qed\)



\subsection*{Merging the label \(0\) squares \(\{A_j\}\) with \(\partial^+_{left}\)}
Using property \((q3)\) above, we merge the label~\(0\) squares \(\{A_j\}_{1 \leq j \leq b}\) iteratively with the cycle~\(\partial^+_{left}.\)



Let \(D_0 := \partial^+_{left}.\) In the first step of the iteration, we consider the label~\(0\) square~\(A_1\) containing the edge \(g_1 \in Q_I.\) The square \(A_1\) lies in the exterior of the cycle \(D_0.\) Using Theorem~\(3\) of Ganesan (2015), we merge \(D_0\) and~\(A_1\) to get a new cycle~\(D_1.\) For \(i \geq 1,\) the cycle \(D_i\) satisfies the following properties.\\\\
\((x1)\) Every label \(1\) square of the plus connected component~\(C^+_{left}\) lies in the interior of the cycle~\(D_i.\) The path \(Q = (h_1,\ldots,h_{n+2}),\) defined in property~\((q1)\) above, is a subpath of \(D_i.\)\\

The path illustrated by the edges \((1,2,\ldots,n+1,n+2)\) in Figure~\ref{wj_fig_plus} is the subpath~\(Q\) and is a subpath of the cycle~\(D_i.\) Intuitively, the second statement of~\((x1)\) true since by property \((q3),\) every square~\(A_j\) we merge lies in the right half plane i.e., to the right of the line~\(x = 0.\) In Figure~\ref{wj_fig_plus}, the line~\(x = 0\) is the infinite line containing the segment \(ak\) of the rectangle~\(R.\) Thus the path~\(Q\) remains unchanged after every iteration.\\\\
\((x2)\) Every edge in \(D_i \setminus Q\) either belongs to the path~\(Q_I\) or belongs to some vacant square~\(A_l, 1 \leq l \leq b.\) Every edge in the path \(Q_I\) either belongs to \(D_i \setminus Q\) or is contained in the interior of the rectangle~\(R.\) For \(1 \leq j \leq b,\) the edge~\(g_{j} \in Q_I\) belongs to~\(D_i \setminus Q\) if and only if the label~\(0\) square~\(A_{j}\) containing~\(g_j\) as an edge lies in the exterior of the cycle~\(D_i.\)\\

The first statement in property~\((x2)\) states that the new cycle obtained,~\(D_i,\) contains only edges of the original path~\(Q_I\) or the edges of the vacant squares in~\(\{A_j\}.\) The second statement in~\((x2)\) states that if an edge of the original path \(Q_I\) does not belong to the cycle~\(D_i,\) then it lies in the interior of \(D_i.\) The final statement in property~\((x2)\) lists a necessary and sufficient condition for an edge of the original subpath~\(Q_I\) to be present in the new cycle~\(D_i.\) This is used to proceed to the next step of the iteration.\\\\
\((x3)\) For~\(1 \leq j \leq i,\) the edge~\(g_j \in Q_I\) and the corresponding label~\(0\) square~\(A_j\) are contained in the interior of the cycle~\(D_i.\) 

Our aim is to merge all the label~\(0\) squares~\(\{A_j\}_{1 \leq j \leq b}\) (see \((q3)\)) iteratively with the cycle~\(\partial^+_{left}.\) Property~\((x3)\) states the progress at the end of the~\(i^{th}\) iteration.


\emph{Proof of~\((x1)-(x3)\) for~\(D_1\)}: From Theorem~\(3\) of Ganesan~(2015), the interior of the new cycle~\(D_1\) contains the interior of the cycle~\(D_0\) and the interior of the square~\(A_1.\) Also the cycle~\(D_1\) consists only of edges of~\(D_0\) and~\(A_1.\) From Theorem~\(2\) of Ganesan~(2015), we have that every label~\(1\) square is contained in the interior of the cycle \(D_0.\) Thus the first statement of~\((x1)\) and the first statement of~\((x2)\) are true.

The second statement of~\((x1)\) is true as follows. We recall that in the proof of the iterative merging algorithm of Theorem~\(3\) of Ganesan (2015), we replace a subpath~\(Z\) of the cycle~\(D_0\) with a subpath~\(Y\) of the square~\(A_1\) so that the following holds. The paths~\(Z\) and~\(Y\) have the same endvertices and the union of the paths \((D_0 \setminus Z) \cup Y\) is the required cycle~\(D_1\) containing the interiors of both~\(D_0\) and~\(A_1.\) This is illustrated in Figure~\ref{wj_fig_plus}, where the edge \(g_1 = ab\) and the label \(0\) square~\(A_1\) is square~\(axyba.\) We replace the subpath~\(Z = ab\) of the path~\(Q_I = abcdefghijk\) with the subpath~\(Y = axyb.\) In the Figure~\ref{wj_fig_plus}, both~\(Y\) and~\(Z\) lie to the right of the line~\(ak\) which is part of the infinite line~\(x=0.\)

By the last statement of property~\((q3),\) every edge in the square~\(A_1\) lies to the right of the line~\(x = 0.\) Since every edge in the path~\(Q\) lies to the left of the line~\(x = 0,\) the cycle~\(D_1\) also contains the subpath~\(Q.\) From the above merging algorithm, we also obtain that every edge in the path~\(Q_I\) either belongs to the cycle~\(D_i\) (in particular, belongs to the path \(D_i \setminus Q\)) or is contained in the interior of~\(D_i.\) This proves the second statement of~\((x2).\)



To prove the last statement of~\((x2)\) suppose that the edge~\(g_j\) of the path~\(Q_I\) belongs to the cycle~\(D_1.\) We recall that~\(g_j\) is common to the label~\(1\) square~\(K_j \in C^+_{left}\) and the label~\(0\) square~\(A_j \in \Lambda^+.\) Therefore either~\(K_{j}\) or~\(A_{j}\) lies in the interior of~\(D_1\) but not both. Since~\(K_j\) lies in the interior of~\(D_1\) (property~\((x1)\)), we have that~\(A_j\) lies in the exterior of~\(D_1.\)

Conversely, suppose that \(g_j \notin D_1.\) Using the second statement of \((x2),\) we therefore have that the edge~\(g_j \in Q_I\) is contained in the interior of~\(D_1.\) Thus both the squares containing~\(g_j\) as an edge belong to the interior of the cycle~\(D_1.\) This proves~\((x2).\) To see~\((x3)\) is true, we have from~\((x2)\) that if~\(g_1 \in D_1,\) then the square~\(A_1\) lies in the exterior of~\(D_1.\) But as mentioned above, the cycle~\(D_1\) contains~\(A_1\) in its interior. Thus~\((x3)\) is true. \(\qed\)


Using property~\((x2),\) we proceed to the next step of the iteration. Fix the least index~\(j\) such that the edge~\(g_{j}\in Q_I\) belongs to the cycle~\(D_1.\) If there is no such~\(j,\) we stop the procedure. If there exists such an index~\(j,\) then from property~\((x2)\) above, the corresponding label~\(0\) square~\(A_{j}\) lies in the exterior of~\(D_1.\) Merge~\(A_{j}\) and~\(D_1\) using Theorem~\(3\) of Ganesan~(2015) to get the new cycle~\(D_{2}.\)


The new cycle~\(D_2\) also satisfies properties~\((x1)-(x3).\)\\
\emph{Proof of~\((x1)-(x3)\) for~\(D_2\)}: It suffices to verify that the squares~\(A_1\) and~\(A_2\) both belong to the interior of the cycle~\(D_2.\) The rest of the proof is as above. Since the cycle~\(D_1\) satisfies \((x1)-(x3)\) we have that~\(A_1\) is contained in the interior of~\(D_1\) and therefore contained in the interior of~\(D_2.\) If~\(A_2\) is also contained in the interior of~\(D_1,\) then we are done.

If not, then the square~\(A_2\) is in the exterior of the cycle~\(D_1\) and the edge \(g_2 \in D_1\) by property \((x2)\) for the cycle~\(D_1.\) In our iteration, we choose the least index~\(j\) such that the edge~\(g_j \in Q_I\) belongs to~\(D_1.\) Thus we choose \(j = 2\) here and we merge~\(A_2\) with~\(D_1\) to get the new cycle~\(D_2.\) This implies that~\(Y_2\) lies in the interior of~\(D_2.\)~\(\qed\)


As before, we use property~\((x2)\) of the cycle~\(D_2\) to proceed to the next step of iteration. This process proceeds for a finite number of steps and the cycle~\(D_i\) obtained at iteration step~\(i \geq 1\) satisfies properties~\((x1)-(x3).\) Let~\(D_{fin}\) denote the final cycle obtained after the procedure stops. The cycle \(D_{fin}\) satisfies the following properties.\\
\((x5)\) The path~\(Q\) (marked by \((1,2,3,\ldots,n+1,n+2)\) in Figure~\ref{wj_fig_plus}) is a subpath of the cycle~\(D_{fin}\) and the subpath~\(Q_{II} = D_{fin}\setminus Q\) has endvertices \((0,0)\) and \((0,n).\)\\
\((x6)\) Every label~\(0\) square in~\(\{A_j\}_{1 \leq j \leq b}\) lies in the interior of the cycle~\(D_{fin}.\)\\

\emph{Proof of \((x5)-(x6)\)}: The property~\((x5)\) is true since at every iteration, the cycle~\(D_i\) satisfies property~\((x1).\) To prove~\((x6),\) we argue as follows. Suppose that \(D_{fin} = D_{i_0}\) is the cycle obtained after \(i_0\) iterations. From the first statement of property~\((x2)\) of the cycle~\(D_{i_0},\) we have that an edge~\(e\) in~\(D_{i_0}\setminus Q\) either belongs to the original path~\(Q_I\) or some label~\(0\) square~\(A_j, 1 \leq j \leq b.\) If~\(e\) belongs to~\(Q_I,\) then~\(e = g_j\) for some~\(1 \leq j \leq b\) and using the final statement of property~\((x2),\) we obtain that the label~\(0\) square~\(A_j\) lies in the exterior of the cycle~\(D_{fin}.\) This means that the iterative procedure would not have terminated. Thus~\(e\) belongs to some label~\(0\) square \(A_l, 1 \leq l \leq b.\)

From the second statement of property~\((x2)\) of the cycle~\(D_{i_0}\) and the discussion in the above paragraph, we have that every edge in the path~\(Q_I\) lies in the interior of~\(D_{i_0}.\) Thus if \(g_j \in Q_I,\) then both the squares containing~\(g_j\) as an edge lie in the interior of the cycle \(D_{i_0}.\) In particular, the label~\(0\) square~\(A_j\) containing~\(g_j\) as an edge lies in the interior of~\(D_{i_0}.\) This proves \((x6).\)\(\qed\)


\subsection*{Constructing the vacant top down crossing using~\(D_{fin}\)}
Using property \((x5),\) we have that the subpath~\(Q_{II}\) intersects the top line~\(y = n\) containing the top edge of the rectangle~\(R\) and it also intersects the bottom line~\(y = 0\) containing the bottom edge of~\(R.\) Let \(Q_{II} = (f_1,\ldots,f_r)\) with the edge~\(f_1\) containing~\((0,n)\) as an endvertex and the edge~\(f_r\) containing~\((0,0)\) as an endvertex.\\\\
\((y1)\) Every edge~\(f_k \in Q_{II}\) is the edge of a unique label~\(0\) square~\(Z_k \in \{A_j\}.\) \\\\
It may happen that multiple edges in \(Q_{II}\) are the edges of the same label~\(0\) square. In any case, the sequence of squares \(\Pi = (Z_1,\ldots,Z_r)\) forms a star connected path since the squares \(Z_k\) and \(Z_{k+1}\) always share the vertex common to \(f_k\) and \(f_{k+1}.\) Moreover, since \(Z_1\) contains \(f_1\) as an edge, the square \(Z_1\) intersects the top line~\(y = n.\) Similarly the square~\(Z_r\) intersects the bottom line~\(y = 0.\)

Let~\(k_1\) be the largest index~\(k \geq 1\) such that the square~\(Z_k\) intersects the top line~\(y = n.\) Similarly let~\(k_2\) be the smallest index~\(k \geq k_1\) such that the square~\(Z_k\) intersects the bottom line~\(y = 0.\) Define the sequence~\(\Pi(k_1,k_2) = (Z_{k_1},Z_{k_1+1},\ldots,Z_{k_2}).\)\\\\
\((y2)\) Suppose there is no occupied plus connected left right crossing of the rectangle~\(R.\) Every square~\(Z_k, k_1 \leq k \leq k_2\) is a vacant square~lying in the interior of the rectangle~\(R.\)

The sequence of squares \(\Pi(k_1,k_2)\) is the desired vacant top down star connected crossing.

\emph{Proof of \((y1)-(y2)\)}:
To prove \((y1),\) suppose there are two label~\(0\) squares~\(A_{j_1}\) and~\(A_{j_2}\) both having the edge~\(f_k.\) The edge~\(f_k \in D_{fin}\) and so one of the squares~\(A_{j_1}\) or~\(A_{j_2}\) necessarily lies in the exterior of~\(D_{fin},\) a contradiction to property~\((x6)\) above.


We prove \((y2)\) in the remaining part of the proof. Suppose that the label~\(0\) square~\(Z_k\) lies in the exterior of the rectangle~\(R\) for some \(k \geq k_1.\) The square \(Z_k \in \{A_j\}\) and every square in~\(\{A_j\}\) lies to the right of the line \(x = 0\) (see property \((q3)\)). Therefore \(Z_k\) cannot lie to the left of the line~\(x=0.\) Suppose that~\(Z_k\) lies above the top line \(y = n.\) Since \(Z_k \in \{A_j\},\) we use property \((q3)\) to obtain that one of the edges of~\(Z_k\) belongs to the original path~\(Q_I.\) From property \((q2),\) we have that every edge in~\(Q_I\) is contained in the rectangle~\(R.\) Therefore the bottom edge~\(b_k\) of~\(Z_k\) belongs to~\(Q_I\) and is contained in the top edge of~\(R.\) The square~\(Z_{k+1}\) shares a corner with the square~\(Z_k.\) If~\(Z_{k+1}\) lies in the interior of the rectangle~\(R,\) then the square~\(Z_{k+1}\) also intersects the top edge of~\(R\) and this contradicts the definition of the index~\(k_1.\)

If \(Z_{k+1}\) lies in the exterior of~\(R,\) then the square~\(Z_{k+1}\) also lies above the top edge of~\(R.\) Arguing as in the previous paragraph, the bottom edge of~\(Z_{k+1}\) is contained in the top edge of~\(R,\) again a contradiction to the definition of the index~\(k_1.\) An analogous argument holds if~\(Z_k\) lies below the bottom edge of the rectangle~\(R.\)

Suppose now that~\(Z_k\) lies to the right of the line~\(x = m\) containing the right edge of~\(R.\) Arguing as above, we have that the left edge~\(l_k\) of~\(Z_k\) belongs to the original path~\(Q_I\) and is contained in the right edge of~\(R.\) From property~\((q3),\) the edge~\(l_k\) is also the edge of a label~\(1\) square~\(Q_k\) belonging to the component~\(C^+_{left}.\) The square~\(Q_k\) is contained in the interior of rectangle~\(R\) and touches the right edge of~\(R,\) a contradiction to property (\ref{lr_c_left_plus}).

From the above, we obtain that each square~\(Z_k, k_1 \leq k \leq k_2\) lies in the interior of the rectangle~\(R\) and has label~\(0.\) Using property (\ref{w_int_plus}), we have that~\(Z_k\) is vacant. This proves~\((y2).\) \(\qed\)








\renewcommand{\theequation}{\thesection.\arabic{equation}}

\subsection*{Acknowledgement}
I thank Professors Rahul Roy and Federico Camia for crucial comments and for my fellowships. I also thank NISER for my fellowship.

\bibliographystyle{plain}

\begin{thebibliography}{10}




\bibitem{boll} B. Bollobas. (2001).
\newblock {\em Modern Graph Theory}.
\newblock {Springer}.

\bibitem{boll} B. Bollobas and O. Riordan. (2006).
\newblock {\em Percolation}.
\newblock {Academic Press}.

\bibitem{gane} G. Ganesan. (2015).
\newblock {Outermost boundaries for star and plus connected components in percolation}.
\newblock {\emph{Arxiv Link}:  http://arxiv.org/abs/1508.06443}.


\bibitem{penrose} M. Penrose. (2003).
\newblock {\em Random Geometric Graphs}.
\newblock {Oxford}.





\end{thebibliography}

\end{document}